\documentclass[12pt,a4paper]{amsart}
\usepackage{amssymb}

\usepackage{graphicx,amssymb,amsfonts,epsfig,amsthm,a4,amsmath,url}
\usepackage[latin1]{inputenc}

\usepackage{amsfonts,amssymb,verbatim}
\usepackage{latexsym}
\usepackage{mathrsfs}
\usepackage{stmaryrd}
\usepackage{xspace}
\usepackage{enumerate, paralist}
\usepackage[usenames,dvipsnames]{color}
 \usepackage{txfonts, pxfonts}

\newtheorem*{thmm}{Theorem}

\newtheorem{thm}{Theorem}[section]
\newtheorem{cor}[thm]{Corollary}
\newtheorem{lem}[thm]{Lemma}
\newtheorem{clai}[thm]{Claim}
\newtheorem{prop}[thm]{Proposition}
\theoremstyle{definition}
\newtheorem{defn}[thm]{Definition}

\newtheorem{conj}[thm]{Conjecture}
\theoremstyle{remark}
\newtheorem{rem}[thm]{Remark}

\numberwithin{equation}{section}

\newcommand{\n}{\mathfrak{n}}
\newcommand{\m}{\mathfrak{m}}
\newcommand{\Z}{\mathbb{Z}}
\newcommand{\HH}{\mathbb{H}}

\newcommand{\E}{\mathbb{E}}

\newcommand{\N}{\mathbb{N}}
\newcommand{\T}{\mathbb{T}}
\newcommand{\R}{\mathbb{R}}

\newcommand{\PP}{\mathcal{P}}

\newcommand{\eps}{\varepsilon}

\newcommand{\var}{\textnormal{var}}
\begin{document}
\title{First passage percolation on nilpotent Cayley graphs}

\author{Itai Benjamini}
\address{The Weizmann Institute, Rehovot, Israel}

\email{itai.benjamini@gmail.com}

\author{Romain Tessera}
\address{Laboratoire de Math\'ematiques, B\^atiment 425\\ Universit\'e Paris-Sud 11, 91405 Orsay, France
}
\email{romain.tessera@math.u-psud.fr}
\date{\today}
\subjclass[2010]{46B85, 20F69, 22D10, 20E22 }
\keywords{First passage percolation, Nilpotent groups, asymptotic cone, invariant random metric on groups}

\baselineskip=16pt

\begin{abstract}
We prove an asymptotic shape theorem for the standard first-passage percolation on Cayley graphs of virtually nilpotent groups. By a theorem of Pansu, the asymptotic cone of a finitely generated nilpotent group is  isometric to a simply connected nilpotent Lie group equipped with some left-invariant Carnot-Caratheodory metric.
Our main result is an extension of Pansu's theorem to  random metrics, where the edges of the Cayley are i.i.d.\ random variable with some finite exponential moment. Based on the companion  work \cite{T}, the proof relies on Talagrand's concentration inequality, together with Pansu's theorem. Adapting an argument from \cite{BKS} we prove a sublinear estimate on the variance for virtually nilpotent groups which are not virtually isomorphic to $\Z$. We further discuss the asymptotic cones of first-passage percolation on general infinite connected graphs: we prove that the asymptotic cones are a.e.\ deterministic if and only the volume growth is subexponential.
\end{abstract}

\maketitle
\tableofcontents


\section{Introduction}

First passage percolation is a model of random perturbation of a given geometry.
In this paper, we shall restrict to the simplest model, where random i.i.d lengths are assigned to the edges of a fixed graph. We refer to  \cite{GK,Ke} for background and references.
A fundamental result (the shape theorem) states that the  random  metric on Euclidean lattices when rescaled by $1/n$, almost surely  converges to a deterministic invariant metric on the Euclidean space \cite{CD,Ke}.
Underlying this theorem is the simple  fact that the graph metric associated to the Euclidean grid when rescaled, converges to the euclidean space equipped with the $\ell^1$-norm. In the world of Cayley graphs, a version of
this last fact holds and characterizes polynomial growth: by a theorem of Gromov \cite{Gr}, groups of polynomial growth are virtually nilpotent, and
by a theorem of  Pansu \cite{Pa}, the rescaled sequence converges in the pointed Gromov-Hausdorff topology to a simply connected nilpotent Lie group equipped with some left-invariant Carnot-Caratheodory metric.
It is therefore natural to ask if when assigning  random i.i.d.\ lengths to a Cayley graph of polynomial growth,
the rescaled metric almost surely converges to a deterministic metric on the Lie group. Establishing this was the original goal of this note. Besides proving it, we also obtain a general statement  on the fluctuations of the distance obtained by first-passage percolation (FPP for short) on general graphs with bounded degree.

Before stating our results, let us describe our general set up.
Consider a connected non-oriented graph $X$, whose set of vertices (resp.\ edges) is denoted by $V$ (resp.\ $E$). We first define the notion of {\it weighted} graph metric on $V$.
For every function $\omega:E\to (0,\infty)$, we equip $V$ with the weighted graph metric $d_{\omega}$, where each edge $e$ has weight $\omega(e)$. In other words, for every $x,y\in V$, $d_{\omega}(x,y)$ is defined as the infimum over all path $p=(e_1,\ldots, e_m)$ joining $x$ to $y$ of $\ell_f(p):=\sum_{i=1}^m\omega(e_i)$. Denote by $d$ the graph metric on $V$,  corresponding to the constant function  $\omega=1$.

Let $\nu$ be a probability measure supported on $[0,\infty)$. The random metric of {\em first passage percolation} consists in choosing the weight $\omega(e)$ independently  according to $\nu$.  Note that $\E d_{\omega}(x,y)$ defines a distance on $V$, that we call the {\it average distance} and denote by $\bar{d}(x,y)$.

A central result in FPP is the following Gaussian concentration inequality due to Talagrand.

\begin{thmm}\cite[Proposition 8.3]{Talagrand}).
Suppose that $\omega(e)$ has an exponential moment: i.e.\ there exists $c>0$ such that $\E \exp(c\omega(e))<\infty$. Then there exists $C_1$ and $C_2$ such that for every graph $X=(V,E)$, for every pair of vertices $x,y$, and for every $u\geq 0$,
\begin{equation}\label{eq:expmoment}
P\left(|d_{\omega}(x,y)-\bar{d}(x,y)|\geq u\right)\leq C_1\exp\left(-C_2\min\left\{\frac{u^2}{d(x,y)},u\right\}\right).
\end{equation}
\end{thmm}

\noindent{\bf Basic assumptions.}
In order to avoid useless repetitions, let us once and for all list the technical assumptions on the edge's length distribution $\nu$, that will be required in most of our statements.

\begin{itemize}
\item {\bf (A1)}
We assume that $\nu$ has an exponential moment, and therefore satisfies (\ref{eq:expmoment}) for some constants $C_1$ and $C_2$ (this assumption can probably relaxed but we choose not to focus on this aspect here). 
\item {\bf  (A2)}
We  also suppose that there exists $a>0$ such that $\bar{d}(x,y)\geq ad(x,y)$ for all $x,y\in V$. 
\end{itemize}

When one works with the standard Cayley graph of $\Z^d$, the second assumption is satisfied exactly when $\nu(\{0\})<p_c$, where $p_c$ is the critical probability of percolation on $\Z^d$ \cite{Ke}.
For more general graphs, we shall also suppose that $\nu(\{0\})<1/k$, where $k$ is an upper bound on the degree of the graph. Indeed,  by \cite[Corollary A2]{T}, this implies condition $(A_2)$. 
Observe that by triangular inequality, $\bar{d}\leq (\E \omega(e))d$.  In the sequel we denote $b:=\E \omega(e).$
It follows that under our second assumption, $d$ and $\bar{d}$ are actually bi-Lipschitz equivalent, more precisely,
\begin{equation}\label{eq:biLip}
ad\leq \bar{d}\leq bd.
\end{equation}

We shall adopt the following notation:  given $v\in V$ and $r>0$, let $\bar{B}(v,r)$ (resp.\ $B_{\omega}(v,r)$) denote the ball of radius $r$ for the average distance $\bar{d}$ (resp.\ for the random distance $d_{\omega}$).

Since this paper addresses to probabilists as well as to geometric group theorists, we start recalling some basic (and less basic) notions of geometric group theory.

\subsection{Cayley graphs and nilpotent groups}
Let $G$ be a finitely generated group, and let $S$ be a finite generating  subset of $G$ such that $S^{-1}=S$. Recall that the Cayley graph $(G,S)$ is defined as follows: the vertex set is $G$ itself, and an edge joins two vertices $g$ and $g'$ if there exists $s\in S$ such that $g'=gs$. We denote by $d_S$ the distance on $G$, obtained by restricting the graph distance to the vertex set $G$. Observe that this distance is left-invariant: if $g,g,k\in G$, then $d_S(kg,kg')=d_S(g,g')$. We shall generally refer to $d_S$ as the word metric associated to $S$.

For group elements $y$ and $z$, let $[y,z]$ denote the commutator element $yzy^{-1}z^{-1}$.  Given two subgroups $A$ and $B$ of the same group $G$, we shall denote by $[A,B]$ the subgroup generated by $[a,b]$ where $a\in A$ and $b\in B$. Let $C^j(G)$ be the descending central series of $G$, i.e.\ let $C^{0}(G)=G$, and $C^{j+1}(G)=[G,C^j(G)]$.  $G$ is $l$-step nilpotent if $C^{l} = \{1\}$ and $C^{l-1} \neq \{1\}$.

Finally, a group is said to be virtually nilpotent if it has a finite index nilpotent subgroup.

\subsection{Nilpotent Lie groups and Carnot-Caratheodory metrics}

We let $N$ be a connected  nilpotent Lie group. Examples of such groups are abelian connected Lie groups such as $\R^d$, but also compact abelian groups such as the $d$-dimensional torus $\T_d\simeq \R^d/\Z^d$. In the sequel we shall only consider {\it simply connected} nilpotent Lie groups, meaning that we exclude the possibility that there is a compact subgroup. This condition is equivalent to requiring that $N$ is homeomorphic to $\R^k$ for some $k$. 
The simplest example of a nilpotent connected simply connected Lie group which is not abelian is the Heisenberg group, whose definition is recalled in the next subsection.

Any connected Lie group $N$ (not necessarily nilpotent) can be endowed with a left-invariant geodesic metric as follows: pick a norm $\|\cdot \|$ on the  tangent space at the neutral element (which identifies with Lie algebra $\n$ of $N$). Now given a smooth path $\gamma$ on $N$, define
the length of $\gamma$
with respect to $\|\cdot \|$  as 
\begin{equation}\label{eq:length}
l(\gamma)=\int_0^1\|\gamma(t)^{-1}\cdot \gamma'(t)\|dt.
\end{equation}
The distance on $N$ is then defined by  
\begin{equation}\label{eq:distance} 
d(x,y)=\inf_{\gamma} \{l(\gamma); \; \gamma(0)=x,\; \gamma(1)=y\}.
\end{equation}
In case $\|\cdot\|$ is euclidean, this actually defines a left-invariant Riemannian metric on $N$.

When $N=\R^n$, the distance defined above is simply the distance induced by the norm $\|\cdot\|$. Observe that an important property of this special case is its ``scale invariance", namely if $x,y\in \R^d$, and if $t>0$, then $d(tx,ty)=\|tx-ty\|=td(x,y)$. 

In order to obtain a suitable generalization of the "scale invariance" on a non-abelian connected simply connected nilpotent Lie group, one needs to  
work with another natural family of left-invariant geodesic distances on $N$, called Carnot-Caratheory metrics. The idea is to start with a norm which is only defined on a subspace of $\n$ and to consider only paths which are tangent to this subspace. More precisely,
let $\n$ be the Lie algebra of $N$, and let $\m$ be a vector subspace supplementing $[\n,\n]$ equipped with a norm $\|\cdot\|$. 
A smooth path $\gamma: [0,1]\to N$ is said to be horizontal if $\gamma(t)^{-1}\cdot \gamma'(t)$ belongs to  $\m$ for all $t\in [0,1]$. The length of $\gamma$
with respect to $\|\cdot \|$ is then defined by (\ref{eq:length}).
It can be shown that since $\m$ generates the Lie algebra $\n$, every pair of points can be joined by a horizontal path (see \cite{Gro}). 
The Carnot-Caratheodory metric associated to $\|\cdot\|$ is defined so that the distance between two points in $N$ is given by (\ref{eq:distance}),
where the infimum is taken over all piecewise horizontal paths $\gamma$. Note that if $N=\R^m$, so that $\m=\n$, then the Carnot-Caratheodory metric is just the distance associated to the norm $\|\cdot\|$. 
We shall see in the next section that in the case of the Heisenberg group, the Carnot Caratheodory metric enjoys some "scale invariance" property, which explains its relevance for the study of limit shape theorems. 

Before stating our main results for nilpotent groups, let us illustrate it in a concrete case. 

\subsection{A limit shape theorem for the Heisenberg group}

Recall that the real Heisenberg group $\HH(\R)$ is defined as the matrix group
$$\HH(\R)=\left\{\left(\begin{array}{cccccc}
1 & u & w \\
0 &  1 & v\\
0 & 0 & 1
\end{array}\right); u,v,w\in \R\right\},$$ and that the discrete Heisenberg  $\HH(\Z)$ sits inside $\HH(\R)$ as the cocompact discrete subgroup
consisting of unipotent matrices with integral coefficients. The group $\HH(\R)$ (resp.\ $\HH(\Z)$) is 2-step nilpotent: indeed, its center, which coincides with its derived subgroup is isomorphic to $\R$ (resp.\ $\Z$), and consists in matrices whose only non-zero coefficient is the top right coefficient. It follows that  $\HH(\R)/[\HH(\R),\HH(\R)]\simeq \R^2$ (and similarly $\HH(\Z)/[\HH(\Z),\HH(\Z)]\simeq \Z^2$).

We equip the group $\HH(\Z)$ with the word metric associated with the finite generating set
$S=\{a^{\pm 1},b^{\pm 1}\}$, where
$$a=\left(\begin{array}{cccccc}
1 & 1 & 0 \\
0 &  1 & 0\\
0 & 0 & 1
\end{array}\right), \; b=\left(\begin{array}{cccccc}
1 & 0 & 0 \\
0 &  1 & 1\\
0 & 0 & 1
\end{array}\right).$$
Consider the one-parameter group $(\delta_t)_{t\in \R^*_+}$ of automorphisms of $\HH(\R)$ defined as follows
$$\delta_t\left(\begin{array}{cccccc}
1 & u & w \\
0 &  1 & v\\
0 & 0 & 1
\end{array}\right)=\left(\begin{array}{cccccc}
1 & tu & t^2w \\
0 &  1 & tv\\
0 & 0 & 1
\end{array}\right).$$
Given a norm $\|\cdot \|$ on $\R^2$, there exists a unique left-invariant Carnot-Caratheodory $d_{cc}$ metric on $\HH(\R)$ that projects to $\|\cdot\|$ and that is scaled by $\delta_t$, i.e.\ such that $d_{cc}(e,\delta_t(g))=td_{cc}(e,g)$ for all $t\in\R^*_+$ and all $g\in \HH(\R)$. Such an automorphism $\delta_t$ (for $t>1$) is called a dilation. In this sense, we can say that $d_{cc}$ is ``scale invariant".

We shall denote by $B_S(g,r)$ the ball of radius $r$ centered at $g\in \HH(\Z)$ for the word metric $d_S$.
Given a left-invariant Carnot-Caratheodory distance $d_{cc}$ on $\HH(\R)$, we let $B_{cc}(g,r)$ denote the ball of radius $r$ centered at $g$ for this distance. Recall that the Hausdorff distance $d_H(A,A')$ between two compact subsets $A$ and $A'$ of $\HH(\R)$ is defined as the minimum over all $r\geq 0$ such that $A\subset [A']_r$ and $A'\subset [A]_r$, where $[A]_r=\{x\in \HH(\R), \; d_{cc}(x,A)\leq r\}$. A sequence $A_n$ is said to Hausdorff converges to $A$ if $d_H(A_n,A)\to 0$. This notion of convergence does not depend on a particular left-invariant Carnot-Caratheodory metric on $\HH(\R)$.

The following theorem is a special case of a theorem of Pansu that we shall recall in complete generality in the next subsection. 
 \begin{thm}\cite{Pa}\label{thm:PansuHei}{\bf (Limit shape for the rescaled discrete Heisenberg group)}
Consider $\HH(\Z)$ equipped with its metric $d_S$. Denote by $d_{cc}$ the Carnot-Caratheodory metric associated to the $\ell^1$-norm on $\R^2.$ Then for every $r>0$, and every $\eps>0$, there exists $n_0$ such that for all $n\geq n_0$,
$$B_{cc}(e,nr(1-\eps))\cap \HH(\Z)\subset  B_S(e,rn)\subset B_{cc}(e,nr(1+\eps)).$$
In particular, $\delta_{1/n}\left(B_S(e,rn)\right)$  Hausdorff converges to  $B_{cc}(e,r)$, as $n\to\infty.$ 
\end{thm}

Our result is  a version of the previous theorem for FPP metrics.

 \begin{thm}\label{thm:Heisen}{\bf (Limit shape for First Passage Percolation on Heisenberg)}
Consider FPP on the Cayley graph $(\HH(\Z),S)$ associated to some measure $\nu$ satisfying both conditions $(A_1)$ and $(A_2)$. Then there exists a (deterministic) Carnot-Caratheodory metric $d_{cc}$ on $\HH(\R)$ such that for every $r>0$, every $\eps>0$,  and a.e.\ every $\omega$, there exists $n_0$ such that for all $n\geq n_0$,
$$B_{cc}(e,nr(1-\eps))\cap \HH(\Z)\subset  B_{\omega}(e,rn)\subset B_{cc}(e,nr(1+\eps)).$$
In particular, $\delta_{1/n}\left(B_{\omega}(e,rn)\right)$ a.e.\ Hausdorff converges to $B_{cc}(e,r)$, as $n\to\infty.$ 
\end{thm}

\subsection{A general result for virtually nilpotent groups}
In order to state a version of Theorems \ref{thm:PansuHei} and \ref{thm:Heisen} for general virtually nilpotent groups, it will be useful to use the notion of Gromov-Hausdorff convergence (see e.g. \cite{BH, BBI, Gromov} for background on this notion).

\begin{defn}\label{GH}
Given a  sequence $X_n$ of compact metric spaces we will say that $X_n$ GH-converges to $X$ if the $X_n$ have uniformly bounded diameter and if there exist maps $\phi_n:X_n\to X$ such that for all $\eps$, then for  $n$ large enough, 
\begin{itemize}
\item every point of $X$ is at $\eps$-distance of a point of $\phi_n(X_n)$;
\item $(1-\eps)d(x,y)-\eps\leq d(\phi_n(x),\phi_n(y))\leq  (1+\eps)d(x,y)+\eps$ for all $x,y\in X_n$.
 \end{itemize}
A sequence of maps $\phi_n$ satisfying these two conditions is called a sequence of GH-approximations of the space $X$.
\end{defn}
GH-convergence naturally extends to (not necessarily compact) locally compact pointed metric spaces (see \cite[Section 3]{Gromov}). 
\begin{defn}  
Given a  sequence $(X_n,o_n)$ of locally compact pointed metric spaces, $(X_n,o_n)$ is said to converge to the locally compact pointed metric space $(X,o)$ if for every $R>0$, the sequence of balls $B(o_n,R)$ GH-converges to $B(o,R)$. 
\end{defn}
A sufficient condition for the sequence  $(X_n,o_n)$ to converge to $(X,o)$ is the existence of a pointed GH-approximations, i.e.\ a sequence of maps $\phi_n:X_n\to X$ such that for all $\eps$, for $n$ large enough
\begin{itemize}
\item $d(\phi_n(o_n),o)\leq \eps$
\item every point of $B(o,1/\eps)$ lies at $\eps$-distance of a point of $\phi_n(B(o_n,1/\eps))$;
\item $(1-\eps)d(x,y)-\eps\leq d(\phi_n(x),\phi_n(y))\leq  (1+\eps)d(x,y)+\eps$ for all $x,y\in (B(o_n,1/\eps))$.
 \end{itemize}

Let us first reformulate Theorem \ref{thm:PansuHei} in this framework: 
consider the sequence of embeddings of $\delta_{1/n}\circ i:\HH(\Z)\to \HH(\R)$, where $i$ is the standard embedding of $\HH(\Z)$ in $\HH(\R)$. This can be interpreted as a sequence of maps $\phi_n$ of pointed metric spaces $(\HH(\Z),d_S/n,e)$ to $(\HH(\R),d_{cc},e)$. Pansu's Theorem can be reformulated by saying that the sequence of maps $\phi_n$ is a sequence of (pointed) GH-approximations. 

In the sequel, we let $G$ be a finitely generated group, $S$ be a finite generating subset, and $(G,S)$ be the corresponding Cayley graph of $G$. 
\begin{thm}\label{thmPrelim:Pansu} \cite{Pa}
Let $G$ be a finitely generated virtually nilpotent group equipped with some finite generating set $S$. Then $(G,d_S/n,1_G)$ converges in the pointed Gromov-Hausdorff topology to some simply connected (Carnot) nilpotent Lie group $N_{\R}$ equipped with some left-invariant Carnot-Caratheodory metric $d_{cc}$.
\end{thm}
Here is our generalization of the previous theorem in the context of FPP.
\begin{thm}\label{thm:Main}{\bf (Asymptotic shape theorem for nilpotent groups)}
Let $G$ be a finitely generated virtually nilpotent group equipped with some finite generating set $S$. Consider FPP on the Cayley graph $(G,S)$ associated to some measure $\nu$ satisfying both conditions $(A_1)$ and $(A_2)$.
Let $(N_{\R},d_{cc})$ be the limit of $(G,d/n,1_G)$ from Theorem \ref {thmPrelim:Pansu}. There exists a left-invariant  Carnot-Caratheodory metric $d'_{cc}$ on $N_{\R}$, which is bi-Lipschitz equivalent to $d_{cc}$ such that for a.e.\ $\omega\in \Omega$,
$(G,d_{\omega}/n,1_G)$ converges in the pointed Gromov-Hausdorff topology to $(N_{\R},d'_{cc})$.
\end{thm}
Giving a formulation of Theorems \ref{thmPrelim:Pansu} and  \ref{thm:Main} in the spirit of Theorem \ref{thm:PansuHei} is possible (see \cite{Pa,Br} for the deterministic case, and it is straightforward to obtain the relevant FPP statement from our arguments). However, this would require a much longer and tedious introduction which we choose to avoid here. The reason for this is that in general, $G$ 
does not sit as a cocompact discrete subgroup inside its rescaled  limit $N_{\R}$. First of all, one would need to pass to a finite index {\it nilpotent and torsion-free} subgroup $H$. Such a group indeed sits as a cocompact discrete subgroup in its {\it Malcev completion} $H_{\R}$, which is a connected simply connected nilpotent Lie group. But even then, it is not always true that $H_{\R}$ is isomorphic to $N_{\R}$. This subtle issue is addressed for instance in \cite{Pa}.

\

The proof of Theorem \ref{thm:Main} goes in two steps: first we use Talagrand's concentration inequality to show that the obvious map $(G,d_{\omega}/n,1_G)\to (G,\bar{d}/n,1_G)$ is almost surely a sequence of Gromov-Hausdorff approximation. This step is completely general: the only geometric property that is used is the fact that the volume of balls in $(G,S)$ grows subexponentially (see Proposition \ref{prop:almostsure}).
The second step consists in showing that $\bar{d}$ is sufficiently close to being geodesic to apply Pansu's theorem  to the sequence $(G,\bar{d}/n,1_G)$. Let us be more specific about this last point. 
\begin{defn}\cite{Pa,Br}
A metric space $X$ is called asymptotically geodesic  (or inner), if for all $\eps>0$, there exists $\alpha$  such that for all $x,y\in X$, there is a sequence $x=x_0,x_1,\ldots, x_m=y$ such that $d(x_i,x_{i+1})\leq \alpha$ for all $0\leq i<m$, and such that 
$$\sum_{i=0}^{m-1}d(x_i,x_{i+1})\leq (1+\eps)d(x,y).$$
\end{defn}
Note that a word metric on a finitely generated group is asymptotically geodesic. 
This notion is motivated by the following strengthening of Theorem \ref{thmPrelim:Pansu}.
Recall that a left-invariant distance on a group $G$ is called proper if for all $r>0$, the ball $B(e,r)$ is finite.
\begin{thm}\label{thmPrelim:PansuInner} \cite{Pa,Br}
Let $G$ be a finitely generated virtually nilpotent group equipped with some left-invariant proper asymptotically geodesic distance $\delta$ on $G$.
 Then $(G,\delta/n,1_G)$ converges in the pointed Gromov-Hausdorff topology to some simply connected (Carnot) nilpotent Lie group $N_{\R}$ equipped with some left-invariant Carnot-Caratheodory metric $d_{cc}$.
\end{thm}

Theorem \ref{thm:Main} now results from the following fact.

\begin{thm}\label{thm:asymptoticgeod}
Let $G$ be a finitely generated virtually nilpotent group equipped with some finite generating set $S$. Consider FPP on the Cayley graph $(G,S)$ associated to some measure $\nu$ satisfying both conditions $(A_1)$ and $(A_2)$. Then the metric space $(G,\bar{d})$ is asymptotically geodesic.
\end{thm}
This is an immediate consequence of the (stronger) \cite[Proposition 1.3]{T}. However for the sake of completeness and since the latter article is not yet published, we reproduce the argument here.
 
\subsection{Asymptotic cone of FPP on graphs with bounded degree}

The first condition to obtain a limit shape theorem is to have relative compactness for the Gromov-Hausdorff topology, which restricts our investigations to graphs with polynomial growth. 
In order to treat more general situations, one needs the notion of asymptotic cone, which is some way of forcing the scaling limit to exist (using some non-principal ultrafilter). These notions are recalled in \S \ref {section:subexp}.
One can then prove a very general result which in some (weak) sense is a far-reaching generalization of the phenomenon observed in Theorem \ref {thm:Main}. 

\begin{thm}\label{thm:ascone}
Let $X=(V,E)$ be a graph with degree at most $k\in \N$, let $o_n$ be a sequence of vertices, $r_n\in \N$ be an increasing sequence, and let $\eta$ be a non-principal ultrafilter. We assume that  $\nu$ is supported on $[a,b]$, with $0<a<b<\infty<$ and that $\nu(\{a\})<1/k$.
Then ``the asymptotic cone is almost surely deterministic", i.e.\ for a.e.\ $\omega$, $$\lim_{\eta}(X,d_{\omega}/r_n,o_n)=\lim_{\eta}(X,\bar{d}/r_n,o_n),$$ if and only if for every $\eps>0$, $$\lim_{\eta}\frac{\log |B(o_n,r_n/\eps)|}{r_n}=0.$$
\end{thm}
Saying that the asymptotic cone is almost surely deterministic amounts to saying that the fluctuations of the metric in the ball of radius $r$ are almost surely ``sublinear", i.e.\ in $o(r)$. For those who do not like ultrafilters and asymptotic cones, we recommend to read the statements of Propositions \ref {thm:ascone} and \ref {prop:exponentialgrowth} which are written in terms of fluctuations.

Theorem \ref {thm:ascone} is the combination of two independent statements: one dealing with the subexponential growth case, and one with the exponential growth case (see Remark \ref{rem}). The first statement (Corollary \ref {cor:SubexpCone}) is a consequence of Talagrand's Theorem, while the second one (Corollary \ref {cor:exponentialgrowth}) is completely elementary. The conclusion of Corollary \ref {cor:exponentialgrowth} is actually stronger than the statement of Theorem \ref {thm:ascone}: roughly speaking it says that the $\omega$-distance in the ball $B(o_n,r_n)$ a.s.\ admits fluctuations of size of the order of $r_n$ about the average distance. We do not know whether this remains true for the distance to the origin.

\subsection{Sublinear upper bound on the variance}
 A straightforward and well-known consequence of Talagrand's theorem is a linear bound on the variance $\var(d_{\omega}(x,y))=O(d(x,y))$ valid for any graph, and sharp for $\Z$ (Kesten first proved it for FPP on $\Z^d$ using martingales \cite{Kes}).
 In \cite{BKS}, the authors manage to improve this linear bound  on $\Z^d$, for $d\geq 2$:
$$\var(d_{\omega}(x,y))\leq C\frac{d(x,y)}{1+\log d(x,y)}.$$
To be more precise, they prove it under the assumption that $\nu(\{a\})=\nu(\{b\})=1/2$. However, in \cite[Theorem 4.4]{BR'}, the same result is proved under much more general assumptions on $\nu$ (including e.g.\ exponential laws). In a subsequent paper, these authors prove a concentration inequality as well \cite[Theorem 5.4]{BR}. All these results rely on the same geometric trick from \cite{BKS}. Although we did not check it, it is likely that they should all be generalized to the setting of Theorem \ref {thm:variance} below. 
\begin{thm}\label{thm:variance}
Assume that $\nu(\{a\})=\nu(\{b\})=1/2$ and consider FPP on some Cayley graph $(G,S)$. Assume that $G$ has a finite index subgroup $G'<G$ whose center $Z(G')$  satisfies the following property:
there exists $\delta>1$ and $c>0$ such that for all $n$
\begin{equation}\label{eq:center}|Z(G')|\cap B_S(e,n)\geq cn^{\delta}.
\end{equation}
Then there exists $C>0$ such that for all $x,y\in G$, one has
\begin{equation}\label{eq:var}\var(d_{\omega}(x,y))\leq C\frac{d(x,y)}{1+\log(1+ d(x,y))}.
\end{equation}
\end{thm}
Let us examine the case of the Heisenberg group: its center is isomorphic to the cyclic subgroup generated by the matrix
$$c=\left(\begin{array}{cccccc}
1 & 0 & 1 \\
0 &  1 & 0\\
0 & 0 & 1
\end{array}\right).$$

Note that $[a^k,b^n]=a^{-k}b^{-n}a^kb^n=c^{nk}$, from which one easily deduces that $|Z(\HH(\Z))|\cap B_S(e,n)|\geq \alpha n^{2}$ for some constant $\alpha>0$. Therefore the previous theorem applies to $\HH(\Z).$ More generally, it is well-known (see e.g.\  \cite{Gui}) that non-virtually abelian nilpotent groups satisfy (\ref {eq:center}) with some $\delta\geq 2$. So Theorem \ref {thm:variance} applies to Cayley graphs of virtually nilpotent groups which are not virtually isomorphic to $\Z$.

\section{Fluctuations of the rescaled distance}
In the sequel, we implicitly assume that we perform FPP on a graph with respect to some measure $\nu$ satisfying $(A_1)$ and $(A_2)$. 
We start by a very simple estimate resulting by union bound from Talagrand's concentration inequality.
\begin{lem}\label{lem:ConvProba}
Let $X=(V,E)$ be a graph, and let $o_n$ be a sequence of vertices.
Let $r_n\geq 0$ be a increasing sequence of integers. 
 For all $\eps>0$,  there exists $n_0>0$  such that for all $n\geq n_0$, the probability that
for all $x,y\in B(o_n,r_n/\eps)$, one has
$$|d_{\omega}(x,y)-\bar{d}(x,y)|\leq \eps r_n,$$
is at least $1-4|B(o_n,r_n/\eps)|^2\exp\left(-\frac{\eps^3 r_n}{64K}\right)$.
\end{lem}

\begin{prop}\label{prop:almostsure} {\bf (Graphs with subexponential growth)}
Let $X=(V,E)$ be a graph, and let $o_n$ be a sequence of vertices. Let $r_n\geq 0$ be an increasing sequence of integers such that for all $\eps>0$,
\begin{equation}\label{eq:subexp}
\log |B(o_n,r_n/\eps)|/r_n=o(r_n)
\end{equation}
There exists a measurable subset of full measure $\Omega'\subset \Omega$ such that
for all $0<\eps\leq \min\{a/2,1/2\}$, and all $\omega\in \Omega'$, there exists $n_0=n_0(\omega,\eps)$ such that for all $n\geq n_0$, for all $x,y\in B(o_n,r_n/\eps)$ with $ d(x,y)\geq \eps r_n$ one has
\begin{equation}\label{eq:uniformfluct}
|d_{\omega}(x,y)-\bar{d}(x,y)|\leq \eps r_n,
\end{equation}
and 
\begin{equation}\label{eq:containment}
B_{\omega}(o_n,r_n)\cup \bar{B}(o_n,r_n)\subset B(o_n, r_n/\eps).
\end{equation}
In particular, the sequence of (tautological) maps $(X,d_{\omega}/r_n,o_n)\to (X,\bar{d}/r_n,o_n)$ is a sequence of pointed GH-approximations.
\end{prop}

\begin{proof}
Note that (\ref{eq:subexp}) is equivalent to the condition that for all $C,c>0$, $$\sum_n  |B(o_n,Cr_n)|^2e^{-cr_n}<\infty.$$
Hence (\ref{eq:uniformfluct}) follows by Borel Cantelli's lemma from Lemma \ref{lem:ConvProba}. Let us check (\ref{eq:containment}). Since $ad\leq \bar{d}$, we have that $\bar{B}(o_n,r_n)\subset B(o_n,r_n/a)\subset B(o_n, r_n/\eps)$. It is therefore enough to show that $B_{\omega}(o_n,r_n)\subset B(o_n,r_n/\eps)$. 
Assume by contradiction that there exists $y\in B_{\omega}(o_n,r_n)$ which does not belong to $B(o_n,r_n/\eps)$. Let $\gamma$ be an $\omega$-geodesic from $o_n$ to $y$, and let $z\in \gamma$ be such that $d(o_n,z)=r_n/\eps$ (assume for simplicity that $r_n/\eps$ is an integer).  On the one hand, we have
$$\bar{d}(o_n,z)\geq ar_n/\eps\geq 2r_n.$$
Which contradicts the following inequality which results from (\ref{eq:uniformfluct}):
$$\bar{d}(o_n,z)\leq d_{\omega}(o_n,z)+\eps r_n<2r_n.$$
Hence we are done.
\end{proof}
We immediately deduce the following statement. 
\begin{cor}\label{cor:limitshape}
Under the assumptions of Proposition \ref{prop:almostsure}, there exists a measurable subset of full measure $\Omega'\subset \Omega$ such that for all $\eps>0$, and all $\omega\in \Omega'$, there exists $n_0=n_0(\omega,\eps)$ such that for all $n\geq n_0$,
$$\bar{B}(o_n,(1-\eps)r_n) \subset B_{\omega}(o_n,r_n)\subset \bar{B}(o_n,(1+\eps)r_n).$$
\end{cor}

We now examine the case of graphs with polynomial growth, for which a more quantitative statement will be needed. 

\begin{prop}\label{cor:fluctuations}{\bf (Graphs with polynomial growth)}
Let $q>0$ and $K>0$, and let  $r_n\in \N$ be an increasing sequence. Then there exists $D>0$, $E>0$, $C>0$, and $n_0\in \N$ such that the following holds.
Let $X=(V,E)$ be a graph and let $o_n$ be a sequence of vertices such that  $|B(o_n,r_n)|\leq Kr_n^q.$ Then 
\begin{itemize}
\item for all $n\geq n_0$,
\begin{equation}\label{eq:uniformFluct}
P\left(\sup_{x,y\in B(o_n,r_n)}|d_{\omega}(x,y)-\bar{d}(x,y)|^2\geq Dr_n\log r_n\right)\leq  Er_n^{-2}.
\end{equation}
\item for a.e.\ $\omega$, there exists $n_1=n_1(\omega)$ such that for $n\geq n_1$,
\begin{equation}\label{eq:uniformVariance}
\sup_{x,y\in B(o_n,r_n)}|d_{\omega}(x,y)-\bar{d}(x,y)|\leq C(r_n\log r_n)^{1/2}.
\end{equation}
\end{itemize}
\end{prop}
\begin{proof}
The second statement follows from the first one via Borel Cantelli's lemma.  Let $D$ be a constant to be determined later. Applying Talagrand's theorem, we obtain that for all large enough $r$, and all $x,y$ such that $d(x,y)\leq r$,
$$P\left(|d_{\omega}(x,y)-\bar{d}(x,y)|^2\geq Dr\log r\right)\leq C_1 \exp\left(-C_2 D\log r\right).$$
Now, letting $D=(2q+2)/(C_2b^2)$, we deduce that for all large enough $r$, and all $x,y$ such that $d(x,y)\leq r$,
$$P\left(|d_{\omega}(x,y)-\bar{d}(x,y)|^2\geq Dr\log r\right)\leq C_1r^{-2q-2}.$$
Hence for $n$ large enough,

$$P\left(\sup_{x,y\in B(o_n,r_n)}|d_{\omega}(x,y)-\bar{d}(x,y)|^2\geq Dr_n\log r_n\right)\leq C_1r_n^{-2q-2}|B(o_n,r_n)|^2\leq C_1Kr_n^{-2}.$$
Hence the first statement follows.
\end{proof}

\section{The average distance is asymptotically geodesic: proof of Theorem \ref{thm:asymptoticgeod}}
This technical section is essentially extracted from \cite{T}. The proofs are repeated in order to make the present paper self-contained.
Throughout this section, we implicitly assume that we perform FPP on a graph with respect to some measure $\nu$ satisfying $(A_1)$ and $(A_2)$.

We start proving the following

\begin{lem}\label{pop:k-geodesic}\cite[Proposition 3.1]{T}
Let $q>0$ and $K>0$, and let  $r_n\in \N$ be an increasing sequence. Then there exists $C>0$ and $n_0$ such that the following holds.
Let $X=(V,E)$ be a graph and let $o_n$ be a sequence of vertices such that  $|B(o_n,r_n)|\leq Kr_n^q.$  Then for all $x,y\in B(o_n,ar_n/(32b))$ and for all $0\leq \lambda\leq 1$, there exists a vertex $z\in B(o_n,r_n)$ such that for $n\geq n_0$,
$$\left|\lambda\bar{d}(x,y)-\bar{d}(x,z)\right|\leq C(r_n\log r_n)^{1/2},$$
and 
$$\left|(1-\lambda)\bar{d}(x,y)-\bar{d}(z,y)\right|\leq C(r_n\log r_n)^{1/2}.$$
\end{lem}

\begin{proof}
By Proposition \ref{cor:fluctuations}, we have
$$P\left(\sup_{x,y\in B(o_n,r_n)}|d_{\omega}(x,y)-\bar{d}(x,y)|^2\geq Dr_n\log r_n\right)\leq Er_n^{-2},$$
with $D=(2d+2)/(C_2b^2)$ (remember that $C_1$ and $C_2$ are the constants appearing in the conclusion of Talagrand's theorem).
Let $n_0$ be the smallest integer so that $Er_n^{-2}<1$. Then for all $n\geq n_0$, there exists $\omega$ (depending on $n$) such that 
\begin{equation}\label{eq:fluctuation}
\sup_{z_1,z_2\in B(o,r_n)}\left|d_{\omega}(z_1,z_2)-\bar{d}(z_1,z_2)\right|\leq D (r_n\log r_n)^{1/2}.
\end{equation}
Assume in addition that $n_0$ is large enough so that $D (r_n\log r_n)^{1/2}\leq ar_n/16$ for all $n\geq n_0$.
Let $\gamma$ be some $\omega$-geodesic between $x$ and $y$. First of all, note  that $\gamma$ cannot escape from the ball $B(o_n,r_n)$. Indeed, suppose there is $1\leq i\leq k$ such that $d(o_n,\gamma(i))= r_n$, then by triangular inequality, $d(x,\gamma(i))\geq r_n/2$, hence $\bar{d}(x,\gamma(i)\geq ar/2$. So (\ref{eq:fluctuation})  implies that 
$$d_{\omega}(x,\gamma(i)) \geq \bar{d}(x,\gamma(i))-D (r_n\log r_n)^{1/2}\geq ar_n/2-ar_n/16\geq ar_n/4.$$ 
which contradicts the fact that  
$$d_{\omega}(x,\gamma(i))\leq d_{\omega}(x,y)\leq \bar{d}(x,y)+ar_n/16\leq ar_n/8.$$ 

By (\ref{eq:fluctuation}), the maximum of $\omega(e)$ over all edges on $\gamma$ is at most $b+D (r_n\log r_n)^{1/2}\leq D'(r_n\log r_n)^{1/2}$ for some $D'>0$. Therefore, one can find a vertex $z$ in $\gamma$ such that 
$$\left|\lambda d_{\omega}(x,y)-d_{\omega}(x,z)\right|\leq D'(r_n\log r_n)^{1/2},$$
and 
$$\left|(1-\lambda)d_{\omega}(x,y)-d_{\omega}(z,y)\right|\leq D' (r_n\log r_n)^{1/2}.$$

But then combining these inequalities with (\ref{eq:fluctuation}), we get 
$$\left|\lambda\bar{d}(x,y)-\bar{d}(x,z)\right|\leq 4D' (r_n\log r_n)^{1/2},$$
and 
$$\left|(1-\lambda)\bar{d}(x,y)-\bar{d}(z,y)\right|\leq 4D'(r_n\log r_n)^{1/2},$$
so that the proposition follows with $c=1/(4D')$.
\end{proof}

We immediately deduce the following corollary.
\begin{cor}\label{cor:SAG*}
Let $q>0$ and $K>0$. Then there exists $C>0$ and $\alpha_0>0$ such that the following holds.
Let $X=(V,E)$ be a graph such that  $|B(o,r)|\leq Kr^q$ for all $o\in V$ and all $r>0$. Then for all $x,y\in V$ such that $r=d(x,y)\geq \alpha_0$, and for all $0\leq \lambda\leq 1$, there exists a vertex $z\in V$ such that
$$\left|\lambda\bar{d}(x,y)-\bar{d}(x,z)\right|\leq C(r\log r)^{1/2},$$
and 
$$\left|(1-\lambda)\bar{d}(x,y)-\bar{d}(z,y)\right|\leq C(r\log r)^{1/2}.$$

\end{cor}

\begin{lem}\label{lem:biadicSG}\cite[Lemma 4.3]{T}
Let $q>0$ and $K>0$. Then there exists $C>0$ and $\alpha_0>0$ such that the following holds.
Let $X=(V,E)$ be a graph such that  $|B(o,r)|\leq Kr^q$ for all $o\in V$ and all $r>0$. 
For all  integer $k\geq 1$, and for all $x,y\in X$ such that $\alpha=\bar{d}(x,y)/2^k\geq \alpha_0$,  there exists  a sequence $x=x_0\ldots, x_{2^k}=y$ satisfying, for all $0\leq i\leq 2^k-1,$
$$\left(1-C\left(\frac{\log \alpha}{\alpha} \right)^{1/2}\right)\alpha \leq \bar{d}(x_i,x_{i+1})\leq \left(1-C\left(\frac{\log \alpha}{\alpha} \right)^{1/2}\right)\alpha,$$
\end{lem}
\begin{proof}We let $r_0=r'_0:=\bar{d}(x,y)$.
We let $n\in \N$ be such that $2^n< r_0\leq 2^{n+1}$.   Assuming that $n$ is large enough  so that $2^{n}\geq \alpha_0$, where $\alpha_0$ is the parameter of Corollary \ref{cor:SAG*}, there exists $z$ such that 
$$r_0/2-C(r_0\log r_0)^{1/2}\leq \bar{d}(x,z),\bar{d}(z,y)\leq r_0/2+C(r_0\log r_0)^{1/2},$$ for some constant $C$. 
We let $r_1=\max\{\bar{d}(x,z),\bar{d}(z,y)\}$ and $r'_1=\min\{\bar{d}(x,z),\bar{d}(z,y)\}$ and apply Corollary \ref{cor:SAG*} to $(x,z)$ and $(z,y)$. Continuing this subdivision process as long as $r'_{k-1}\geq \beta$, we find a sequence $r_1,\ldots, r_k,\ldots,$ satisfying 
\begin{equation}\label{eq:rk}
r_{k}\leq \frac{1}{2}\left(r_{k-1}+C(r_{k-1}\log r_{k-1})^{1/2}\right),
\end{equation}
and 
\begin{equation}\label{eq:rk}
r_{k}'\geq \frac{1}{2}\left(r'_{k-1}-C(r'_{k-1}\log r'_{k-1})^{1/2}\right),
\end{equation}
and a sequence of finite sequences of vertices $x=z_0(k),\ldots,z_{2^{k}}(k)=y$ such that $$r'_k\leq d(z_i(k),z_{i+1}(k))\leq r_{k},$$
for all $0\leq i<2^{k}-1.$

\begin{clai}\label{lem:biadic}
There exists a constant $A$ such that for all $k$ such that $r'_k\geq \beta$, $$A^{-1}2^{-k}r_0\leq r'_k\leq r_k\leq A2^{-k}r_0.$$
\end{clai}
\begin{proof}
Let us first prove the right inequality, the other one being similar. Let $k\geq 2$, and observe that $$r_k\leq \frac{1}{2}\left(r_{k-1}+Cr_{k-1}^{2/3}\right).$$  We do the following change of variable: $A_k=2^{-k}r_k$ (note that $A_k\geq 1$).  We have $$A_k\leq A_{k-1}+CA_{k-1}^{2/3}2^{-k/3}\leq A_{k-1}\left(1+C2^{-k/3}\right),$$ 
from which we easily deduce that $A_k$ is bounded by some $A$ only depending on $C$.
\end{proof}
In what follows, we assume that $k$ is such that $r'_k\geq \beta$.
We deduce from the lemma and from the fact that $r_k\geq r_0/2^k\geq 2^{n-k}$ (which follows by triangular inequality) that 
\begin{eqnarray*}
r_{k} & \leq & \frac{1}{2}\left(r_{k-1}+C(r_{k-1}\log r_{k-1})^{1/2}\right)\\
         & \leq & \frac{r_{k-1}}{2}\left(1+C\left(\frac{\log r_{k-1}}{r_{k-1}}\right)^{1/2}\right)\\
         & \leq & \frac{r_{k-1}}{2}\left(1+C(n-k+A+1)^{1/2}2^{-(n-k)/2}\right)\\
         & \leq & \frac{r_{k-1}}{2}\left(1+C'(n-k)^{1/2}2^{-(n-k)/2}\right)
\end{eqnarray*}
for some constant $C'$.  
Taking the log and using that $\log(1+x)\leq x$ for $x\geq 0$, we have 
\begin{eqnarray*}
\log (2^kr_{k}/r_0) & \leq & C'\sum_{i=1}^k(n-i)^{1/2}2^{-(n-i)/2}\\
                             & \leq & C'\sum_{j\geq n-k}j^{1/2}2^{-j/2}\\
                             & \leq & C"(n-k)^{1/2}2^{-(n-k)/2},
\end{eqnarray*}
for some constant $C">0$. Remember that $r_k'$, and therefore $A2^{n-k}$ is supposed to be larger than $\alpha_0$. Up to enlarging $\alpha_0$ if necessary we can assume that  $C"(n-k)^{1/2}2^{-(n-k)/2}\leq 1$. Then, using that $\exp(x)\leq 1+2x$, for all $0\leq x\leq 1$, we deduce that there exists a constant $C$ such that
\begin{equation}\label{eq:rkabove}
r_k\leq 2^{-k}\left(1+2C"(n-k)^{1/2}2^{-(n-k)/2}\right)r_0\leq \alpha\left(1+C\left(\frac{\log \alpha}{\alpha} \right)^{1/2}\right),
\end{equation}
where $\alpha=r_0/2^k$.
We prove similarly that 
\begin{equation}\label{eq:rk'below}
r'_k\geq \alpha\left(1-C\left(\frac{\log \alpha}{\alpha} \right)^{1/2}\right).
\end{equation}
We let $x=x_0=z_0(k),\ldots, x_{2^k}=z_{2^k}(k)=y$. We deduce from (\ref{eq:rkabove}) and (\ref{eq:rk'below}) that there exists a constant $C$ such that for every $0\leq i\leq 2^k-1$,
$$\left(1-C\left(\frac{\log \alpha}{\alpha} \right)^{1/2}\right)\alpha\leq d(x_i,x_{i+1})\leq \left(1-C\left(\frac{\log \alpha}{\alpha} \right)^{1/2}\right)\alpha,$$
where $\alpha=\bar{d}(x,y)/2^k$. So Lemma \ref{lem:biadicSG} follows. 
\end{proof}

\begin{proof}[Proof of Theorem \ref{thm:asymptoticgeod}] Gven $0<\eps\leq 1$, choose $\alpha$ large enough so that $$C\left(\frac{\log \alpha}{\alpha} \right)^{1/2}\leq \eps.$$ Then it follows from Lemma \ref{lem:biadicSG} that there exists a sequence $x=x_0\ldots, x_{2^k}=y$ satisfying
\begin{itemize}
\item $d(x_i,x_{i+1})\leq 2\alpha$, for all $0\leq i\leq 2^k-1$; 
\item  $
\sum_{i=0}^{2^k-1} d(x_i,x_{i+1})\leq \bar{d}(x,y)(1+\eps)$.
\end{itemize}
So we are done.
\end{proof}

\section{Asymptotic shape theorem for FPP on nilpotent groups}
We give the concluding steps of the proof of Theorem \ref{thm:Main}. As explained in the introduction, it will result from Proposition \ref{prop:almostsure}, Theorem \ref{thmPrelim:PansuInner} and Theorem \ref{thm:asymptoticgeod}. 

\begin{proof}[Proof of Theorem \ref{thm:Main}]
By Proposition  \ref{prop:almostsure}, for a.e.\ $\omega$, the sequence of tautological maps $(G,d_{\omega}/n,e)\to (G,\bar{d}/n,e)$ is a sequence of GH-approximation. Hence it is enough to show that 
there exists a left-invariant  Carnot-Caratheodory metric $d'_{cc}$ on $N_{\R}$, which is bi-Lipschitz equivalent to $d_{cc}$ such that $(G,\bar{d}/n,1_G)$ converges in the pointed Gromov-Hausdorff topology to $(N_{\R},d'_{cc})$.
Note that by $(A_1)$, the distances $\bar{d}$ and $d$ are bi-Lipschitz equivalent, so that if it exists, $d'_{cc}$ is automatically  bi-Lipschitz equivalent to $d_{cc}$. On the other hand, by Theorem \ref{thm:asymptoticgeod}, $(G,\bar{d})$ is asymptotically geodesic, so Theorem \ref{thmPrelim:PansuInner} ensures the existence of $d'_{cc}$.
\end{proof}

We now turn to the proof of Theorem \ref{thm:Heisen}. 
Let us first recall its deterministic counterpart which is due to Pansu (see also \cite{Br}).

 \begin{thm}\cite{Pa}\label{thm:PansuHeiGen}
Consider $\HH(\Z)$ equipped with a left-invariant, proper, asymptotically geodesic metric $\delta$. Then there exists a  Carnot-Caratheodory $d_{cc}$ on $\HH(\R)$ such that  for every $r>0$, and every $\eps>0$, there exists $n_0$ such that for all $n\geq n_0$,
$$B_{cc}(e,nr(1-\eps))\cap \HH(\Z)\subset  B_{\delta}(e,rn)\subset B_{cc}(e,nr(1+\eps)).$$
\end{thm}

\begin{proof}[Proof of Theorem \ref{thm:Heisen}]
We first deduce from Theorem \ref{thm:asymptoticgeod} that Theorem \ref{thm:PansuHeiGen} applies to the average distance $\bar{d}$. Theorem \ref{thm:Heisen} therefore follows from Corollary \ref{cor:limitshape}.
\end{proof}

\section{Asymptotic cones of FPP on graphs with bounded degree}\label{section:subexp}

\subsection{Ultralimits, asymptotic cone, and Gromov-Hausdorff convergence}\label{section:HG}
\label{subsection:prelim1}

 First recall that an ultrafilter (see \cite{Com}) is a map from $\eta:\PP(\N)\to \{0,1\}$, such that $\eta(\N)=1$, and which is ``additive" in the sense that $\eta(A\cup B)=\eta(A)+\eta(B)$ for all $A$ and $B$ disjoint subsets of $\N$. 
Ultrafilters are used to ``force" convergence of bounded sequences of real numbers. Namely, given such a sequence $a_n$, its limit is the only real number $a$ such that for every $\eps>0$ the subset $A$ of $\N$ of integers $n$ such that $|a_n-a|<\eps$ satisfies $\eta(A)=1$. In this case, we denote $\lim_{\eta} a_n=a$. An ultrafilter is called non-principal if it vanishes on finite subsets of $\N$. Non-principal ultrafilters are known to exist but this requires the axiom of choice. In the sequel, let us fix some non-principal ultrafilter $\eta$. 

\begin{defn}
Given a sequence of pointed metric spaces $(X_n,o_n)$, its ultralimit with respect to $\eta$ is the quotient of 
$$\{(x_n)\in \Pi_n X_n, \; \exists C>0,\; \forall n, \; d(x_n,o_n)\leq C\}$$
 by the equivalence relation $x_n\sim y_n$ if $\lim_{\eta}d(x_n,y_n)= 0$. It is equipped with a distance defined by $d((x_n),(y_n))=\lim_{\eta}d(x_n,y_n).$
 \end{defn}
It is a basic fact that a sequence $a_n\in \R$ converging  to $a$ satisfies $\lim_{\eta}a_n=a$. This fact actually extends to ultralimits of metric spaces:

\begin{lem}\label{lemPrelim:ultra}\cite[Exercice 5.52]{BH} If a sequence of pointed metric spaces converges in the pointed GH sense to $X$, then its ultralimit with respect to $\eta$ is isometric to $X$. 
\end{lem}

In this sense, ultralimits generalize the notion of (pointed) GH-limits. To see why the latter is much more restrictive, recall  the following 
\begin{prop}\label{lemPrelim:compact}
(Gromov's compactness criterion, \cite[Theorem 5.41]{BH}) A sequence of compact metric spaces $(X_n)$ is relatively Gromov Hausdorff compact if and only if the $X_n$'s have bounded diameter, and are ``equi-relatively compact": for every $\eps>0$, there exists $N\in \N$ such that for all $n\in \N$, $X_n$ can be covered by at most $N$ balls of radius $\eps.$
\end{prop}

Let us close this short discussion with the notion of asymptotic cone, which formalizes the idea of ``scaling ultralimit" of a metric space.
\begin{defn}
Given a metric space $X$, a sequence of points $o_n\in X$, an increasing sequence $r_n$ going to $\infty$,  and a non-principal ultrafilter $\eta$, the asymptotic cone of $X$ relative to this data  is the ultralimit $\lim_{\eta}(X,d/r_n,o_n)$.
\end{defn}

\subsection{Graphs with subexponential growth: ``fluctuations vanish in the asymptotic cone".}

Let us mention an immediate consequence of Proposition \ref{prop:almostsure}. 
\begin{cor}\label{cor:SubexpCone}
Under the assumption of Proposition \ref{prop:almostsure}, there exists a measurable subset $\Omega'$ of full measure such that for all $\omega\in \Omega'$ and all non-principal ultrafilter $\eta$,
$$``\lim_{\eta}(X,d_{\omega}/r_n,o_n)=\lim_{\eta}(X,\bar{d}/r_n,o_n)",$$
in the sense that, for all $x_n,y_n\in V $ such that $d(o_n,x_n)= O(r_n)$ and $d(o_n,y_n)= O(r_n)$, $$\lim_{\eta}d_{\omega}(x_n,y_n)/r_n=\lim_{\eta}\bar{d}(x_n,y_n)/r_n.$$ In other words, ``the asymptotic cone is almost surely deterministic".
\end{cor}
As a special case of the previous corollary, we deduce that the asymptotic cone of  FPP on a Cayley graph with subexponential growth is almost surely deterministic.

It is important to make a clear distinction between the strong statement of Corollary \ref{cor:SubexpCone}, and the following much weaker one, which is true on any graph.
\begin{prop}\label{prop:sequencewise}
 Let $X=(V,E)$ be any graph, and for every $n\in \N$, let $o_n\in V$ and $r_n\geq 0$, and let $\eta$ be a non-principal ultrafilter. Then for all $x_n,y_n\in V $ such that $d(o_n,x_n)= O(r_n)$ and $d(o_n,y_n)= O(r_n)$, there exists a measurable subset of  full measure $\Omega'$ (depending on the sequence) such that for all $\omega\in \Omega'$,
 $$\lim_{\eta} d_{\omega}(x_n,y_n)/r_n=\lim_{\eta}\bar{d}(x_n,y_n)/r_n.$$
 \end{prop}
\begin{proof}
This is a consequence of the Lemma \ref{lem:ConvProba}.
\end{proof}


 \subsection{Graphs with exponential growth: ``fluctuations remain non-trivial in the asymptotic cone".}
In this subsection, we shall make the assumption that $\nu$ is supported on an interval $[a,b]$ with $0<a<b<\infty$. 

Note that for first-passage percolation on the r-regular tree for $r\geq 3$, it is easy to see that the asymptotic cone of FPP is not deterministic. We can use the fact that the random distance between two vertices $x,y$ in the tree is only determined by the edges along the unique geodesic between them: this distance is therefore the sum of $n:=d(x,y)$ independent random variables. The average distance $\bar{d}(x,y)$ is   equal to $cd(x,y)$, where $c\in (a,b)$ is the expected length of a given edge.  The probability that $d_{\omega}(x,y)$ is less than --say $(a+c)d(x,y)/2$ (resp.\ more than $(c+b)d(x,y)/2$) decays (at most) exponentially with $n$. On the other hand, there are at least exponentially many pairs of disjoint geodesics of length $n$ in a ball of radius $kn$, for $k\geq 2$. Moreover, the exponential exponent can be made as large as we want by increasing $k$: for instance, for all $x$ in the sphere of radius $(k-1)n$, pick a geodesic joining $x$ to a point of the sphere of radius $kn$. It follows for a.e.\ $\omega$, one can find in the asymptotic cone a pair of distinct points whose $\omega$-distance is strictly less (or strictly larger) than the average distance\footnote{The same argument adapts to non-elementary hyperbolic graphs. To generalize the previous argument, one uses the fact that  there exists $C>0$ such that for all $\omega$ and for every pair of points $x,y$, there is a geodesic (say for the word metric) $\gamma$ between $x$ and $y$ whose $C$-neighborhood  contains any $d_{\omega}$-geodesic between $x$ and $y$.  To conclude that there exist fluctuations  of linear size (both above and below the average distance), one needs to produce exponentially many ``independent" pairs of points at distance $n$ in a ball of radius $\lesssim n$: this follows for instance by considering a quasi-isometrically embedded $3$-regular tree.}.

For general graphs of exponential growth (even Cayley graphs), we do not know whether it is possible to exhibit fluctuations {\it above} the average distance in the asymptotic cone. However, it is  possible to show that it always has fluctuations {\it below} the average distance: More precisely, the following proposition says that if the growth is exponential, then a.s.\  one can find in the asymptotic cone  pairs of distinct points whose $\omega$-distance are ``as close as possible to the minimal possible distance $a(dx,y)$". Provided that the average distance is bounded away from this minimal distance (see Lemma \ref {lem:average}), this implies that  FPP admits ``random fluctuations of linear size", which are therefore visible in the asymptotic cone. 

\begin{prop}\label{prop:exponentialgrowth}
Let $X=(V,E)$ be a (not necessarily connected) graph with bounded degree, let $o_n$ be a sequence of vertices. Assume that there exists an increasing sequence $r_n\in \N$ such that  $\log |B(o_n,r_n)|\geq cr_n$ for some constant $c>0$. Then there exists a measurable subset $\Omega"$ of full measure with the following properties.
For all $\omega\in \Omega"$, for all $\eps>0$, there exists $\tau>0$ and $x_n,y_n\in V$ such that $d(o_n,x_n)= O(r_n)$ and $d(o_n,y_n)= O(r_n)$  and such that for all $n$ large enough,
$$\tau r_n\leq  d_{\omega}(x_n,y_n)\leq (a+\eps)d(x_n,y_n).$$
Moreover, if $\nu(\{a\})>0$, then one can take $\eps=0$.
\end{prop}

Before proving this proposition, let us restate it in terms of asymptotic cones.
\begin{cor}\label{cor:exponentialgrowth}
Let $X=(V,E)$ be a (not necessarily connected) graph with  bounded degree, let $o_n$ be a sequence of vertices and let $\eta$ be a non-principal ultrafilter. Assume that there exists an increasing sequence $r_n\in \N$ such that $\log |B(o_n,r_n)|\geq cr_n$ for some constant $c>0$.  Then there exists a measurable subset $\Omega"$ of full measure with the following properties.
For all $\omega\in \Omega"$, for all $\eps>0$, there exist $x_n,y_n\in V$ such that $d(o_n,x_n)= O(r_n)$ and $d(o_n,y_n)= O(r_n)$  such that
$$0<\lim_{\eta} d_{\omega}(x_n,y_n)/r_n\leq (a+\eps)\lim_{\eta} d(x_n,y_n)/r_n.$$
Moreover, if $\nu(\{a\})>0$, then one can take $\eps=0$.
\end{cor}

\begin{proof}
Note that since the degree of $X$ is bounded, there exists $C$ such that
\begin{equation}\label{eq:expGrowth}e^{cr_n}\leq |B(o_n,r_n)|\leq e^{Cr_n}.
\end{equation}
Let $\lambda=c/2C$, so that $|B(o_n,\lambda r_n|\leq e^{cr_n/2}$.
We now consider a subset $X_n$ of $B(o_n,r_n)$ whose points are pairwise at distance at least $(c/4C)r_n$ apart and which is maximal for this property.
It follows that $$B(o_n,r_n) \subset \bigcup_{x\in X_n}B(x,(c/2C)r_n),$$
from which we deduce that
$$|B(o_n,r_n)|\leq |X_n|e^{cr_n/2}.$$
Thus we deduce that
$$|X_n|\geq e^{cr_n/2}.$$
 We let $0<\lambda<1$ to be determined later and let $k_n=[\lambda (c/8C)r_n]$. Observe that the balls $B(x,k_n)$, for $x\in X_n$ are pairwise disjoint. So one can pick for every $x\in X_n$ a point $y_x$ at distance $k_n$ from $x$, and a geodesic (for $d$) $\gamma_x$ between them.
The probability that all edges of $\gamma_x$ have $\omega$-length at most $(1+\eps)a$ is at least $\nu([a,a(1+\eps)])^{k_n}$. Since the paths $\gamma_x$ are disjoint, these events are independent, so that the probability that one of them has $\omega$-length at most $(a+\eps)k_n$ is at least $$1-(1-\nu([a,a(1+\eps)])^{k_n})^{|X_n|}\geq 1-(1-\nu([a,a(1+\eps)])^{\lambda (c/8C)r_n})^{\exp(cr_n/2)}.$$
Recall that given two sequences such that  $u_n\to 0$ and $v_n\to \infty$, one has $(1-u_n)^{v_n}\leq \exp(-u_nv_n)$. On the other hand, by taking $\lambda$ small enough (depending on $\eps$, unless $\nu(\{a\})>0$), one can ensure that $e^{cr_n/2}\nu([a,a(1+\eps)])^{\lambda (c/8C)r_n}\geq \exp(c'r_n)$ for some $c'>0$.  Therefore for this choice of $\lambda$, the above probability tends to $1$ as $n$ tends to infinity very quickly (in particular the probability of the complement event is summable). This is enough to ensure the existence of a measurable subset of full measure $\Omega"$ such that for all $\omega\in \Omega"$, there is a sequence $x_n\in X_n$ such that for $n$ large enough, $d_{\omega}(y_{x_n},x_n)\leq (a+\eps) k_n'=a(1+\eps) d(y_{x_n},x_n)$.
This proves the first part of the proposition with $y_n=y_{x_n}$.
\end{proof}

To finish the proof of Theorem \ref {thm:ascone}, we need the following lemma.

\begin{lem}\label{lem:average}\cite[Lemma 2.1]{T}
Let $X=(V,E)$ be a graph of degree $\leq q$. Assume that $\nu(\{a\})<1/q$. Then there exists $c>a$ such that $\bar{d}(x,y)\geq a'd(x,y)$ for all $x,y\in V$.
\end{lem}

\begin{rem}\label{rem}
To conclude the proof of Theorem \ref {thm:ascone}, let us remark that in the proof of Corollary  \ref {cor:SubexpCone} (resp.\ in Corollary \ref {cor:exponentialgrowth}) the condition $\log |B(o_n,r_n)|=o(r_n)$ (resp.\ $\log  |B(o_n,r_n)|\geq cr_n$) only needs to hold $\eta$-almost surely.
\end{rem}

\section{Upper bound on the variance}

The proof of Theorem \ref {thm:variance} is a simple generalization of the proof of \cite[Theorem 1]{BKS} (which deals with the case of $\Z^d$, $d\geq 2$). We shall sketch its proof, following the same order as in \cite{BKS}, but only providing justifications when the argument needs to be adapted to our more general setting.
To simplify the exposition, we shall assume that $\delta\geq 2$: for $\delta>1$, the idea is the same but the details are slightly more tedious. Moreover, in the case of nilpotent groups we are interested in, one can always assume $\delta\geq 2$ as recalled in the introduction. In this section, we will denote $1_G$ for the neutral element of $G$, keeping the letter $e$ for the edges. Remember, since this will play a crucial role in this proof that the graph structure on $(G,S)$ is defined by saying that two elements (i.e.\ vertices) $g$ and $g'$ and joined by an edge if there exists $s\in S$ such that $g'=gs^{\pm 1}$. Hence, the action by left-translations of $G$ on itself preserves the graph structure and thus the metric.

Following \cite{BKS}, let us fix $g\in G$, and consider the random variable
$f(\omega):=|g|_{\omega},$ where $|g|_{\omega}$ denotes the $\omega$-distance from the neutral element to $g$. We shall also denote $|g|=d(1_{G},g)$, where $d$ is the word metric on $G$.
For every $\omega$, we pick some $\omega$-geodesic $\gamma$ from $1_G$ to $g$. For every $\omega$ and every edge $e\in E$ we denote $\sigma_e\omega$ the configuration which is different from $\omega$ only in the $e$-coordinate.
We start remarking that
\begin{equation}\label{eq:gamma}
\sum_{e\in E}P(e\in \gamma)\leq (b/a)|g|.
\end{equation}

We then fix $m=[d(1_G,g)^{1/4}]$ and consider the function $g_m:\{a,b\}^{m^2}\to \{1,\ldots, m\}$ constructed in \cite[Lemma 3]{BKS}. Let $\Sigma=\{0,1\}\times \{1,\ldots, m^2\}$, $ \tilde{ \Omega}:=\{a,b\}^{\Sigma}$. We define an injective map  $\psi:\tilde{ \Omega} \to  \{1,\ldots, m\}^2$ by
$$\psi(x)=\left(g_m(x_{0,1},\ldots,x_{0,m^2}),g_m(x_{1,1},\ldots,x_{1,m^2})\right).$$

We let $C\geq 1$ be large enough so that $|Z(G')\cap B_S(1_G,Cm)|\geq m^2$, and we pick some injective map from $j:\{1,\ldots, m\}^2\to Z(G')\cap B_S(1_G,Cm)$.
Let $z:=j\circ\psi$.
We can now define $\tilde{f}$ as a map from $\{a,b\}^{\Sigma}\times \{a,b\}^{E}$ to $Z(G')\cap B_S(1_G,Cm)$ by
$$\tilde{f}(x,\omega)=d_{\omega}(z(x),z(x)g).$$
The first important estimate from \cite{BKS} is
\begin{equation}\label{eq:var(tilde)}
\var(f)\leq \var(\tilde{f})+O\left(m\sqrt{\var{\tilde{f}}}\right)+O(m^2).
\end{equation}
If $z$ commutes with $g$, as
$d(1_G,z)=d(g,gz)\leq Cm$, we deduce by triangular inequality that $|f-\tilde{f}|\leq 2bCm$, which implies (\ref {eq:var(tilde)}). More generally, one needs  that $gz=z'g,$ for some $z'\in Z(G')$ such that $d_S(e,z')=O(m)$. This is guaranteed by the following lemma, after noticing that up to replacing $G'$ with the intersection of all its images by automorphisms of $G$, we can assume that $G'$ is a characteristic subgroup of $G$, whose center is therefore normal in $G$: hence $z'=gzg^{-1}\in Z(G')$.
\begin{lem}
Assume $G'$ is characteristic. There exists some constant $C$, such that for all $g\in G$ and $z\in Z(G')$,
$d_S(e,gzg^{-1})\leq Cd_S(e,z)$.
\end{lem}
\begin{proof}
Note that the action by conjugation of $G$ on $Z(G')$ factors through $G/G'$ which is finite. Let $F\subset G$ be a set of representatives of $G/G'$, and let $C=\max_{g\in F,s\in S}d_S(e,gsg^{-1})$. Let $z\in Z(G')$ of length $n$, and let $z=s_1\ldots,s_n$, where $s_i\in S$. Given $g\in G$, there exists $h\in F$ such that $gzg^{-1}=hzh^{-1}$. Thus we have
$$g^{-1}zg=(hs_1h^{-1})\ldots (hs_nh^{-1}),$$
so the lemma follows by triangular inequality.
\end{proof}

Define $$I_e(\tilde{f}):=P\left(\tilde{f}(x,\omega)\neq \tilde{f}(x,\sigma_e(\omega)\right)=2P\left(\tilde{f}(x,\omega)< \tilde{f}(x,\sigma_e(\omega)\right).$$
Then one needs to show that
\begin{equation}\label{eq:Ie}
I_e(\tilde{f})=O(d(1_G,g)^{-1/4}),
\end{equation}
and
\begin{equation}\label{eq:Ie'}
\sum_e I_e(\tilde{f})=O(d(1_G,g)).
\end{equation}
The rest of the proof is identical to \cite{BKS} so we will not repeat it.
Note that if the pair $(x,\omega)\in \tilde{\Omega}$ satisfies $\tilde{f}(x,\sigma_e \omega)> \tilde{f}(x,\omega)$, then $e$ must belong to every geodesic between $z$ ($=z(x)$) and $zg$. Hence conditioning on $z$ and translating both $\omega$ and $e$ by $z^{-1}$ gives 
\begin{equation}\label{eq:IP}
I_e(\tilde{f})\leq 2P\left(z^{-1}e\in \gamma\right).
\end{equation}
Let $Q$ be the set of edges $e'$ such that $P(z^{-1}e=e')>0$. Note that $e'$ lies in the $B(1_G,Cm)$-orbit of $e$, so that once again the lemma ensures that $Q$ has diameter in $O(m)$. It results that $\gamma\cap Q$ contains $O(m)$ edges.
We now need the following property of $g_m$ (\cite[Lemma 3]{BKS}):  $$\max_{y}P(g_m(x)=y)=O(1/m),$$
from which we deduce that
$$\max_{z_0}P(z=z_0)=O(1/m^2).$$
Conditioning on $\gamma$ and summing over the edges in $\gamma\cap Q$, we get
$$P\left(e\in z\gamma|\gamma\right)= O(1/m).$$
Consequently (\ref{eq:IP}) and the choice of $m$ give (\ref {eq:Ie}).
Also, (\ref {eq:gamma}) implies
$$\sum_{e\in E}P\left(z^{-1}e\in \gamma|z\right)\leq (b/a)|g|.$$
Combining this with (\ref {eq:IP}) gives (\ref {eq:Ie'}), so we are done.

\section{Remarks and questions}

\subsection{More general distributions}
It would be interesting to investigate whether our results survive to  non-trivial correlations between edges lengths. Note that in some sense, Talagrand's exponential concentration estimate is far too strong for Theorem \ref{thm:Main}: actually a polynomial decay with a large exponent would be enough to beat the (polynomial) growth rate of the group. This suggests that one should be able to use weaker estimates possibly allowing some weak correlations.

For groups, one can consider a different type of generalization:  given an ergodic $G$-probability space $(\Omega,P)$, an invariant random metric (IRM) on $G$ is a measurable map $G\times G\times \Omega\to \R_+$, $(g,h,\omega)\to d_{\omega}(g,h)$, such that for a.e.\ $\omega\in \Omega$, $d_{\omega}(\cdot,\cdot)$ is a distance on $G$, and that satisfies the equivariance condition: for a.e.\ $\omega$, and all $g,h_1,h_2\in G$,
$$d_{g\omega}(gh_1,gh_2)=d_{\omega}(h_1,h_2).$$ Clearly FPP is a special case of IRM, where the space $\Omega$ is $[a,b]^E$ equipped with the product probability. Observe that in this case, the action of $G$ on $\Omega$, induced by its action of $E$, is ergodic (actually even mixing).

One may wonder under what conditions on an IRM is the asymptotic cone of $(G,d_{\omega},e)$ almost surely deterministic.
In the special case of virtually nilpotent groups, one may ask whether $(G,d_{\omega},e)$ converges in the pointed Gromov-Hausdorff topology to a connected Lie group equipped with an invariant Carnot-Caratheodory metric. Classical proofs of the limit shape theorem for $\Z^d$ are based on the subadditive ergodic theorem, which allows to treat very general IRM (see \cite{Bj}  for the most general known statement). Unfortunately, we were not able to exploit the subadditive ergodic theorem for non-virtually abelian nilpotent groups: this only gives us that distances along certain ``horizontal" directions are asymptotically deterministic, but for instance in the case of Heisenberg, it is not clear under what conditions distances in the direction of the center do not have large fluctuations.

Let us discuss this in more details.
Recall that the proof of Theorem \ref{thm:Main} splits into two independent parts: one consists in proving a concentration phenomenon, namely that the identity map $(G,d_{\omega}/n,e)\to(G,\bar{d}/n,e)$ induces a sequence of Gromov-Hausdorff approximations (recall that $\bar{d}=\E d_{\omega}$).  This might remain true under very general assumptions on $d_{\omega}$, and in particular it  may not require $d_{\omega}$ to be geodesic, not even asymptotically. This contrasts with the second step,  consisting in proving that $(G,\bar{d}/n,e)$ converges, which does require $\bar{d}$ to be asymptotically geodesic: indeed,   conversely, if $(G,\bar{d}/n,e)$ converges to some geodesic metric space, then $\bar{d}$ must be asymptotically geodesic. 
On the other hand one can exhibit invariant metrics on the Heisenberg group which are not asymptotically geodesic and yet quasi-isometric to the word metric. Moreover such a metric $d$ can be chosen so that $(G,d/n,e)$ does not converge at all \cite[Remark A.6.]{C}.

\subsection{Sublinear variance}
The proof of the sublinear estimate on the variance (Theorem \ref{thm:variance}) uses the fact that the group has a large center.  By contrast, we know that for $\Z$, or more generally on a tree, the variance grows linearly (this can easily be extended to Gromov-hyperbolic graphs).  We suspect that --at least in the context of Cayley graph-- the fact that the variance is sublinear might be related to the fact that no asymptotic cone has cut points (a cut point has the property that when we remove it, the space becomes disconnected). We propose the following more modest conjecture
\begin{conj}
Suppose $G$ is the direct product of two infinite finitely generated groups, then (\ref {eq:var}) is satisfied for all Cayley graphs of $G$.
\end{conj}
A particularly interesting case is the direct product of the $3$-regular tree $T$ with $\Z$: in this case, \cite{BM} have managed to prove that $\E(|d_{\omega}-\bar{d}|)$ is tight in the $\Z$-direction. There is some reason to believe that in the $T$-direction the variance should behave as for $\Z^2$ (since geodesics are likely to remain at bounded distance from the direct product of a geodesic in $T$ times $\Z$). Overall, the variance should be even smaller for $T\times \Z$ than for $\Z^2$, where it is classically conjectured to be of the order of $n^{2/3}$ (we refer to \cite{BKS} and \cite{GK} for a more detailed discussion concerning $\Z^2$). Another interesting example is the product of two $3$-regular trees, for which no sublinear estimate is known at the moment.


\subsection{RWRE on virtually nilpotent Cayley graphs}

The FPP shape theorem and the rate of convergence are a statements regarding
large scale  metric homogenization of local  random metric perturbations.  
Similarly to the path we took here for FPP,  
it is of interest to consider the random walk, heat kernel and Green functions homogenization in the context of virtually nilpotent Cayley graphs. Extending the work 
from lattices in Euclidean spaces, studied in PDE under the name of homogenization and in probability theory under the name RWRE (random walk in random environment).

\bigskip
\footnotesize

\end{document}